\documentclass[a4paper]{amsart}

\usepackage{indentfirst}
\usepackage{a4wide}
\usepackage{amssymb}
\usepackage{amsmath}
\usepackage{graphics}
\usepackage{epsfig}

\newcommand{\RR}{\mathbb R}

\newcommand{\TT}{\mathbb T}
\newfont{\fung}{cmff10 at 15pt}
\newcommand{\ee}{\varepsilon}
\newcommand{\dv}{\hbox{ div}}

\newcommand{\la}{\Lambda}
\newcommand{\laa}{\Lambda^{\alpha}}
\newcommand{\las}{\Lambda^{s}}

\newcommand{\ep}{\epsilon}
\newcommand{\ff}{\mathbb}
\newcommand{\pa}{\partial}

\catcode`@=11 \@addtoreset{equation}{section}
\renewcommand\theequation{\thesection.\@arabic\c@equation}
\catcode`@=12

\newtheorem{theorem}{Theorem}[section]
\newtheorem{lemma}[theorem]{Lemma}
\newtheorem{proposition}[theorem]{Proposition}
\newtheorem{corollary}[theorem]{Corollary}
\newtheorem{remark}[theorem]{Remark}
\newtheorem{definition}[theorem]{Definition}

\def\eop{{\ \vrule height 3pt width 3pt depth 0pt}}

\def\cu{{\rm curl}}

\def\sig{\hbox{\rm sign}}

\date{April 4, 2008}

\begin{document}

\title[Porous media with fractional diffusion]{Incompressible flow in porous media with fractional diffusion}

\author[A. Castro, D. C\'ordoba, F. Gancedo and R. Orive]{\'Angel Castro, Diego C\'ordoba, Francisco Gancedo and Rafael Orive}

\address{Angel Castro and Diego C\'ordoba \hfill\break \indent
Departamento de Matem\'aticas \hfill\break\indent Instituto de
Ciencias Matem\'aticas \hfill\break\indent Consejo Superior de
Investigaciones Cient\'{\i}ficas \hfill\break\indent Serrano 123,
28006 Madrid, Spain.}\email{ {\tt
angel.castro@imaff.cfmac.csic.es} {\sc and} {\tt
dcg@imaff.cfmac.csic.es}
  }

\address{Francisco Gancedo \hfill\break \indent
Department of Mathematics \hfill\break\indent University of
Chicago \hfill\break\indent 5734 University Avenue, Chicago, IL
60637, USA .}\email{{\tt fgancedo@math.uchicago.edu}}

\address{Rafael Orive \hfill\break \indent
Departamento de Matem\'aticas \hfill\break\indent Facultad de
Ciencias \hfill\break\indent Universidad Aut\'onoma de Madrid
\hfill\break\indent Crta. Colmenar Viejo km.~15,  28049 Madrid,
Spain.}\email{{\tt rafael.orive@uam.es}}

\thanks{The authors
were partially supported by the grant {\sc MTM2005-05980}
of the {\sc MEC (Spain)} and S-0505/ESP/0158 of the CAM (Spain).
The fourth author was partially supported by the grant {\sc MTM2005-00714} of the {\sc MEC (Spain)}.}

\keywords{Flows in porous media
\\
\indent 2000 {\it Mathematics Subject Classification.} 76S05, 76B03, 65N06.}

\begin{abstract}
In this paper we study the heat transfer with a general
fractional diffusion term  of an incompressible
fluid in a porous medium governed by Darcy's law.  We show formation of singularities
with infinite energy and  for finite energy we obtain
existence and uniqueness results of strong solutions for
the sub-critical and critical cases. We prove global existence of weak solutions for different cases.  Moreover, we obtain the
decay of the solution in $L^p$, for any $p\geq2$,
and the asymptotic behavior is shown. Finally, we prove the existence of an attractor in a weak sense and, for the sub-critical dissipative case with $\alpha\in ( 1,2]$, we obtain the existence of the global attractor for the solutions in the space $H^s$ for any $s > (N/2)+1-\alpha$.
\end{abstract}

\maketitle

\section{Introduction}

We use Darcy's law to model the flow velocities, which yields the following relationship between the liquid discharge
(flux per unit area) $v\in \RR^N$ and  the pressure $$v=-k\left(\nabla p+g\gamma T\right),
$$
where  $k$ is the matrix  medium permeabilities in the different directions respectively
divided by the viscosity, $T$ is the liquid temperature, $g$ is the acceleration due to gravity and the
vector $\gamma\in \RR^N$ is the last canonical vector $e_N$.
While the Navier-Stokes and the Stokes systems are both microscopic equations,
Darcy's law yields a macroscopic description of a flow in the porous medium \cite{bear}. To simplify the notation, we consider $k=g=1$.

In this paper we study the transfer of the heat with a general diffusion term in an incompressible flow. The system which we consider is the following (for more details  see \cite{n-b}):
\begin{eqnarray}\label{mpv}
& &\displaystyle\frac{\partial T}{\partial
t}+v\cdot\nabla T=-\nu\laa T,\\
& &\label{darcy}
v=-\left(\nabla p+\gamma T\right),\\
& & \label{incom}
\dv v=0,
\end{eqnarray}
where $\nu>0$, and the operator $\la^\alpha$ is given by
$\la^\alpha\equiv (-\Delta)^{\alpha/2}$. We will treat the cases $0\leq \alpha \leq 2$ and denote it by DPM. The
case $\alpha = 1$ is called the critical case, the case $1<\alpha \leq 2$ is
sub-critical and the
case $0\leq \alpha < 1$ is super-critical. Roughly speaking, the critical and super-critical cases are
mathematically harder to deal with than the sub-critical case.

In \cite{fab1} P. Fabrie investigates a system of partial differential equations
describing the natural convection in a porous medium under a gradient of temperature, which is obtained by coupling the energy equation and the
Darcy-Forchheimer equation. He proved existence, uniqueness and regularity of the evolution problem as well as the existence of stationary solutions
for the two-dimensional case. Moreover, a regularity theorem is established and a uniform estimate in time of the second-order space derivatives of
the solutions of the three-dimensional case is given. In \cite{fab-nico}  the authors consider the large-time behavior of solutions to the system
\begin{eqnarray*}
&&\gamma v\sb t+v+\nabla p-Ra\sp *\gamma T=0,\\&& {\rm div} v=0,\\ &&T\sb t-\Delta T+v·\nabla T-v\sb3=0,
\end{eqnarray*}
describing the natural convection in a porous medium filling a bounded domain in $\RR^3$.
The asymptotic behavior of the solutions was studied using the concept of an exponential attractor, i.e.
a compact finite-dimensional set invariant under the flow
associated with the system, and uniformly exponentially attracting all the
trajectories within a bounded absorbing set. The main results include the existence of exponential attractors as well as their strong continuity in
a singular (adiabatic) limit $\gamma\to 0$.  In \cite{ly-titi}, using a different method (Galerkin), it was established a global existence and uniqueness
result for the strong solutions of the three-dimensional B\'enard convection problem in a porous medium. Furthermore, a Gevrey class regularity is
obtained for the finite-dimensional attractor of the system. Later, in \cite{oli-titi}, the authors deduce the $H^1\times H^2$ regularity of the attractor.
Combining this with a Fourier splitting method, they were able
to establish the real analyticity of solutions in the attractor.

More recently the Boussinesq approximation of the equations of
coupled heat and fluid flow in a porous medium is studied in
\cite{e-f-z}. This system corresponds to
\eqref{mpv}--\eqref{incom} with $\alpha=2$. They showed that the
corresponding system of partial differential equations possesses a
global attractor. They give lower and upper bounds of the
Hausdorff dimension of the attractor depending on a physical
parameter of the system, namely the Rayleigh number of the flow.

Next, we rewrite the system \eqref{mpv} to obtain the velocity in terms of $T$. The 2D inviscid case is shown in \cite{cgo}. Due to the incompressibility condition, we have that $
\Delta v=- \cu(\cu v)$. Then by computing the curl of the curl of Darcy's law (\ref{darcy}), we get
$$
\Delta v=\left(\frac{\partial^2 T}{\partial x_1\partial x_3}, \frac{\partial^2 T}{\partial x_2\partial x_3},
-\frac{\partial^2 T}{\partial x_1^2}-\frac{\partial^2 T}{\partial x_2^2}\right)
$$
Taking the inverse of the Laplacian
$$
v=\frac{1}{4\pi}\int\limits \frac{1}{|x-y|}\left(\frac{\partial^2 T}{\partial x_1\partial x_3},
\frac{\partial^2 T}{\partial x_2\partial x_3},-\frac{\partial^2 T}{\partial x_1^2}-\frac{\partial^2 T}{\partial x_2^2}\right)dy
$$
and integrating by parts we obtain
\begin{eqnarray}
\label{2.4} v(x,t)= -\frac23 (0,0,T(x,t))+\frac1{4\pi}PV\int_{\RR^3}K(x-y)T(y,t)dy,\qquad x\in \RR^3,
\end{eqnarray}
where
$$
K(x)=\left(3\frac{x_1x_3}{|x|^5},3\frac{x_2x_3}{|x|^5},\frac{2x_3^2-x_1^2-x_2^2}{|x|^5}\right).
$$

\smallskip

In sections 2 through 4 we consider the case where the spatial
domain  can be either the whole $\RR^N$ or the torus $\TT^N$ with
periodic boundary condition.

In section 2 we obtain results of existence  of strong solutions
of the system  \eqref{mpv}--\eqref{incom} under the hypothesis of
regular initial data $T_0\in H^s$ with $s>0$ and $\alpha\in
(1,2]$. The case $\alpha=2$ was studied also in \cite{e-f-z}. For
the supercritical case $\alpha\in[0,1)$, there is global existence
for small initial data $T_0\in H^s$ with $s> N/2+1$. Also, in the
critical case $\alpha=1$, the existence of strong solutions is
obtained as in \cite{c-v,k-n-v} for the critical dissipative
quasi-geostrophic equation.

In section 3  we present results of global existence of weak
solutions. We prove a generalization of the classical
Leray-Prodi-Serrin condition for the uniqueness of the solutions
and obtain global existence and uniqueness for the subcritical
case. In section~4, we obtain the decay of
the solutions of \eqref{mpv}--\eqref{incom} in $\RR^N$ and $\TT^N$
for the $L^p$-norms.

Since we are dealing with a dissipative system,  we study in
section~5 some attracting properties of the solutions of
\eqref{mpv}--\eqref{incom} in $\TT^N$ with a source term $f$ time
independent:
 \begin{eqnarray}\label{mpv2}
& &\displaystyle\frac{\partial T}{\partial
t}+v\cdot\nabla T+\nu\laa T=f.
\end{eqnarray}
It is easy to see that $\overline{T}$, the mean of the solution $T$ of
(\ref{mpv2}), satisfies
$$
\frac d{dt}\overline{T} =\overline{f}
$$
where
$$
\overline{T}(t) =\int\limits_{\TT^N}T(x,t) dx,\quad\hbox{ and}\quad
\overline{f} =\int\limits_{\TT^N}f(x) dx.
$$
Therefore, without loss of generality, we can assume that $f$ and
$T$ are always mean zero. In particular, we prove the existence of
a global attractor in the set of the weak solutions  with the
weak-topology of $L^2_{loc}(0,\infty, L^2(\TT^N))$ and a global
classic compact attractor, connected and maximal in $H^s$ with
$s>N/2+1-\alpha$ in the topology of the strong solutions.

In  section~6 we present results of local existence and blow up of  solutions with infinite energy for $\nu= 0$ and
for $\nu > 0$ in the case $\alpha = 1,2$.

\smallskip\

\section{Strong solutions}\label{s4}

Here we show global existence results of the DPM system
\eqref{mpv}--\eqref{incom} in the sub-critical case. We use a
maximum principle for the $L^p$ norm of the solutions of DPM,
\begin{equation}\label{4.1}\| T(t)\|_{L^p}\leq\| T_0\|_{L^p}\qquad\hbox{with }\:1\leq p\leq \infty,\end{equation}
which is a consequence of $\nabla \cdot v=0$ and the following positivity lemma (see \cite{Resnick} and \cite{cor2}):

\begin{lemma}
For $f,\laa f\in L^p$ with $0\leq\alpha\leq 2$ and $1\leq p$, it is satisfied
\begin{equation}\label{4.2}
\int |f|^{p-1}\sig(f)\laa f\,dx\geq 0.
\end{equation}
\end{lemma}

\smallskip

\begin{theorem}\label{sub}
Let $ T_0\in H^s\cap L^p$ with $s>0$ and $N/(\alpha-1)<p<\infty$. Then, there exists  $T\in C([0,\infty);H^s)$, solution of DPM with $1<\alpha\leq2$.
\end{theorem}

\noindent {\bf Proof.} For a solution of DPM we have
$$
 T_t=-\dv(v T)-\nu\laa T.
$$
We use the equality $\partial_{x_i}=\la(R_i)$, where $R_i$ is the Riesz transforms (see \cite{St1}), to get
$$
\frac{1}{2}\frac{d}{dt}\|\las T\|_{L^2}^2=-\int\la^{s+\frac{\alpha}{2}} T\la^{s+1-\frac{\alpha}{2}}
(R_i(v_i T))\,dx-\nu\|\la^{s+\frac{\alpha}{2}} T\|_{L^2}^2.
$$
 H\"older inequality and the Calderon-Zygmund inequalities for the Riesz transforms (see \cite{St1}) give
\begin{align*}
\frac12\frac{d}{dt}\|\las T\|_{L^2}^2&\leq\|\la^{s+\frac{\alpha}{2}} T\|_{L^2}
\|\la^{s+1-\frac{\alpha}{2}}(v T)\|_{L^2}-\nu\|\la^{s+\frac{\alpha}{2}} T\|_{L^2}^2.
\end{align*}
By the estimate for the operator $\las$ applied to the product of functions (see \cite{taylor}) for $s>0$
\begin{equation}\label{resnisk}
\|\las(fg)\|_{L^r}\leq C(\|f\|_{L^{q'}}\|\las g\|_{L^q}\|+\|g\|_{L^{q'}}\|\las f\|_{L^q})\qquad 1<r<q'\leq\infty,\quad
\frac{1}{r}=\frac{1}{q}+\frac{1}{q'},
\end{equation}
we have for $q'=p$,
$$
\frac12\frac{d}{dt}\|\las T\|_{L^2}^2 \leq
C\|\la^{s+\frac{\alpha}{2}} T\|_{L^2}(\|v\|_{L^p}\|\la^{s+1-\frac{\alpha}{2}}
 T\|_{L^q}+\| T\|_{L^p}\|\la^{s+1-\frac{\alpha}{2}} v\|_{L^q}) -\nu\|\la^{s+\frac{\alpha}{2}} T\|_{L^2}^2,
$$
with $(1/p)+(1/q)=(1/2)$. Now, since $v$ satisfies \eqref{2.4}, again we apply the Calderon-Zygmund inequalites obtaining
$$
\|v\|_{L^p}\leq C\| T\|_{L^p},\qquad\|\la^{s+1-\frac{\alpha}{2}} v\|_{L^q}\leq C\|\la^{s+1-\frac{\alpha}{2}}
 T\|_{L^q},
$$
and \eqref{4.1} gives
$$
\frac12\frac{d}{dt}\|\las T\|_{L^2}^2\leq C\| T_0\|_{L^p}\|\la^{s+\frac{\alpha}{2}} T\|_{L^2}
\|\la^{s+1-\frac{\alpha}{2}}  T\|_{L^q}-\nu\|\la^{s+\frac{\alpha}{2}} T\|_{L^2}^2.
$$
The inequality for the Riesz potential (see \cite{St1})
\begin{equation}\label{rieszpotential}
\|I_{\beta}(f)\|_{L^q}\leq C\|f\|_{r},\quad 0<\beta<N,\quad 1<r<q<\infty,\quad \frac1q=\frac1r-\frac\beta N,\quad I_{\beta}=\la^{-\beta},
\end{equation}
for $r=2$ and $\beta=N/p$, yields
$$\|\la^{s+1-\frac{\alpha}{2}}  T\|_{L^q}\leq C\|\la^{s+1-\frac{\alpha}{2}+\beta}  T\|_{L^2}.$$
We take $p>N/(\alpha-1)$, and therefore $1+\beta< \alpha$, so that
$$
\|\la^{s+1-\frac{\alpha}{2}+\beta}  T\|_{L^2}\leq
\|\la^{s+\frac{\alpha}{2}} T\|_{L^2}^{\gamma}\|\la^{s} T\|_{L^2}^{1-\gamma},
$$
with $\gamma=(2-\alpha+2\beta)/\alpha<1$.
Applying the last two inequalities we obtain
$$
\frac12\frac{d}{dt}\|\las T\|_{L^2}^2\leq
C\| T_0\|_{L^q}\|\la^{s+\frac{\alpha}{2}} T\|_{L^2}^{1+\gamma}\|\la^{s} T\|_{L^2}^{1-\gamma}
-\nu\|\la^{s+\frac{\alpha}{2}} T\|_{L^2}^2,
$$
and Young's inequality gives
$$
\frac12\frac{d}{dt}\|\las T\|_{L^2}^2
\leq C(\nu,\| T_0\|_{L^q})\|\la^{s} T\|_{L^2}^2.
$$
Furthermore, we have
$$
\|\las T\|_{L^2}(t)\leq\|\las T_0\|_{L^2}e^{Ct}.
$$
From this a priori inequality together with the energy estimates argument we can conclude the global existence result.\quad\eop

\smallskip

\begin{theorem}
 Let $0\leq\alpha\leq 1$ be given and assume that $T_0\in H^s$, $s>(N-\alpha)/2+1$. Then there is a time $\tau=\tau(\|\la^sT_0\|)$ so that there exists a unique solution to DPM with $T\in C([0,\tau),H^s)$.
\end{theorem}

\noindent {\bf Proof.} Since the fluid is incompressible we have
for $s>(N-\alpha)/2+1$
\begin{align*}
\frac{1}{2}\frac{d}{dt}\|\las T\|^2_{L^2} &=-\int\las T\las(v\nabla T)\,dx-\nu\|\la^{s+\frac\alpha2}T\|^2_{L^2}\\
&=-\int\las T(\las(v\nabla T)-v\las(\nabla T))\,dx-\nu\|\la^{s+\frac\alpha2}T\|^2_{L^2}\\
&\leq  C\|\las T\|_{L^2}\|\las(v\nabla T)-v\las(\nabla T)\|_{L^2}-\nu\|\la^{s+\frac\alpha2}T\|^2_{L^2}.
\end{align*}
Using  the following estimate (see \cite{Kato-Ponce})
$$
\|\las(fg)-f\las(g)\|_{L^p}\leq C \left(\|\nabla f\|_{L^\infty}\|\la^{s-1} g\|_{L^p}+ \|\las f\|_{L^p}
\|g\|_{L^\infty}\right)\qquad1<p<\infty,
$$
we obtain for $p=2$
\begin{equation*}
\frac{1}{2}\frac{d}{dt}\|\las T\|^2_{L^2}\leq C(\|\nabla v\|_{L^\infty}+ \|\nabla T\|_{L^\infty})
\|\las T\|_{L^2}^2-\nu\|\la^{s+\frac\alpha2}T\|^2_{L^2}.
\end{equation*}
Applying Sobolev estimates we get
\begin{equation*}
\frac{1}{2}\frac{d}{dt}\|\las T\|^2_{L^2}\leq C(\|T_0\|_{L^2}+ \|\la^{N/2+1+\ee}T\|_{L^2})\|\las T\|_{L^2}^2
-\nu\|\la^{s+\frac\alpha2}T\|^2_{L^2}
\end{equation*}
and taking $\ee=s+\frac\alpha2-\frac N2-1$ it follows
\begin{equation}\label{lesuper}
\frac{1}{2}\frac{d}{dt}\|\las T\|^2_{L^2}\leq C(\frac{1}{\nu}+1)(\|T_0\|^2_{L^2}+ \|\las T\|^2_{L^2})\|\las T\|_{L^2}^2.
\end{equation}
Local existence is a consequence of the above a priori inequality.

Let us consider two solution $T^1$ and $T^2$ of DPM
with velocity $v^1$ and $v^2$ respectively, and equal to the
initial datum $T^1(x,0)=T^2(x,0)=T_0(x)$. If we denote $T=T^1-T^2$
and $v=v^1-v^2$, we have
\begin{equation*}
\frac{1}{2}\frac{d}{dt}\|T\|^2_{L^2}\leq -\int T v\cdot \nabla T^1\,dx-\nu\|\la^{\frac\alpha2}T\|^2_{L^2}.
\end{equation*}
For $\alpha=0$ Calderon Zygmund and Sobolev estimates give
\begin{equation*}
\frac{1}{2}\frac{d}{dt}\|T\|^2_{L^2}\leq C\|T\|^2_{L^2}\|\nabla T^1\|_{L^\infty}\leq C\|T\|^2_{L^2}\|T^1\|_{H^s}.
\end{equation*}
Inequality \eqref{lesuper} implies that $\|T^1\|_{H^s}(t)$ is
locally bounded. Furthermore, we can conclude
$$
\|T\|^2_{L^2}(t)\leq \|T\|^2_{L^2}(0)\exp(C\int_0^t\|T^1\|_{H^s}(\sigma)d\sigma),
$$
which yields  uniqueness.

The case $\alpha>0$ is treated  differently, we have
\begin{equation*}
\frac{1}{2}\frac{d}{dt}\|T\|^2_{L^2}\leq \|T\|_{L^2}\|T\|_{L^q}\|\nabla T^1\|_{L^p}-\nu\|\la^{\frac\alpha2}T\|^2_{L^2},
\end{equation*}
with $q=2N/(N-\alpha)$, and $p=2N/\alpha$. Since
\eqref{rieszpotential} we obtain
$$
\|T\|_{L^q}\leq C\|\la^{\frac\alpha2}T\|_{L^2},\quad \|\nabla T^1\|_{L^p}\leq C\|\la^{1+\frac N2-\frac\alpha2}T^1\|_{L^2}\leq C\|T^1\|_{H^s}
$$
and finally
\begin{equation*}
\frac{1}{2}\frac{d}{dt}\|T\|^2_{L^2}\leq \frac C\nu \|T\|^2_{L^2}\|T^1\|_{H^s}^2.\quad\eop
\end{equation*}

\begin{remark}
For the supercritical cases ($0\leq\alpha< 1$), we have the same
criterion as \cite{cgo} for the formation of singularities in
finite time. In fact, we have that
 $ T\in
C([0,\tau]; H^s)$ with $s>N/2+1$ for any $\tau>0$ if, and only if,
$$
\int_0^\tau\|\nabla T\|_{BMO}(t)\,dt<\infty.
$$
\end{remark}
For small initial data we obtain the following global existence result for the supercritical case.

\begin{theorem} Let $\nu>0$, $0\leq\alpha< 1$, and the initial datum satisfies the
smallness assumption
$$\| T_0\|_{H^s}\leq\frac{\nu}{C},\qquad s>N/2+1,$$ for $C$ a fixed constant. Then, there exists a unique
solution of \eqref{mpv}--\eqref{incom} in $C([0,\infty);H^s).$
\end{theorem}

\noindent {\bf Proof.} We multiply \eqref{mpv} by $\la^{2s}T$  and, by the Sobolev embedding, we get
\begin{align*}
\frac{1}{2}\frac{d}{dt}\|\las T\|^2_{L^2}&\leq C(\|\nabla v\|_{L^\infty}+ \|\nabla  T\|_{L^\infty})\|\las
 T\|_{L^2}^2-\nu\|\la^{s+\frac{\alpha}{2}} T\|_{L^2}^2\\
&\leq C(\| T\|_{L^2}+\|\las T\|_{L^2})\|\las  T\|_{L^2}^2 -\nu\|\la^{s+\frac{\alpha}{2}} T\|_{L^2}^2
\end{align*}
with $s>2$. Thus, we have
\begin{align*}
\frac{1}{2}\frac{d}{dt}(\| T\|^2_{L^2}+\|\las T\|^2_{L^2})&\leq
-\nu\|\la^{\frac{\alpha}{2}} T\|_{L^2}^2+C(\| T\|_{L^2}+\|\las T\|_{L^2})\|\las  T\|_{L^2}^2
-\nu\|\la^{s+\frac{\alpha}{2}} T\|_{L^2}^2.
\end{align*} Since
$$
\|\las T\|^2_{L^2}\leq\|\la^{\frac{\alpha}{2}} T\|^2_{L^2}+ \|\la^{s+\frac{\alpha}{2}} T\|^2_{L^2},
$$ we obtain
$$
\frac{1}{2}\frac{d}{dt}(\| T\|^2_{L^2}+\|\las T\|^2_{L^2}) \leq
\|\las T\|_{L^2}^2(C(\| T\|_{L^2}+\|\las T\|_{L^2}) -\nu)\leq0
$$
by  the assumption of the smallness of the initial datum. \quad\eop


\smallskip\

In the critical case $\alpha=1$, we state the following regularity result.

\begin{theorem}\label{caf-vas}
Let $T$ be a solution to the system \eqref{mpv}--\eqref{incom}. Then $T$ verifies the level set energy inequalities, i.e., for every $\lambda>0$
$$
\int\limits_{\RR^N} T_\lambda^2(t_2,x) dx + \int^{t_2}_{t_1}\int\limits_{\RR^N} |\Lambda^{1/2}T_\lambda|^2 dx dt\leq \int\limits_{\RR^N} T_\lambda^2(t_1,x) dx, \qquad 0<t_1<t_2
$$
where $T_\lambda=(T-\lambda)_+$. It yields that for every $t_0>0$ there exists $\gamma>0$ such that $T$ is bounded in $C^\gamma([t_0,\infty)\times\RR^N)$.
\end{theorem}

With this result we show that the solutions of the diffusive
porous medium with initial $L^2$  data and critical diffusion
$(-\Delta)^{1/2}$ are locally smooth for any space dimension. The
proof is analogous to the critical dissipative quasi-geostrophic
equation that it is shown in \cite{c-v}. Analogous result  can be
obtained using the ideas of \cite{k-n-v} to show that  the
solutions in 2D with periodic $C^\infty$ data remain $C^\infty$
for all time.

\smallskip\

\section{Weak solutions}

In this section we prove the global existence of weak solutions for
DPM with $0< \alpha\leq 2$. First we give the definition of weak
solution.
\begin{definition}
The active escalar $T(x,t)$ is a weak solution of DPM if for any $\varphi\in
C^{\infty}_c([0,\tau]\times\RR)$ with $\varphi(x,\tau)=0$, it follows:
\begin{equation}\label{sd}
0=\int\limits_{\RR^N}T_0(x)\varphi(x,0)dx+  \int\limits_0^\tau\int\limits_{\RR^N}T(x,t)
\left(\partial_t\varphi(x,t)+ v(x,t)\cdot\nabla\varphi(x,t)-\nu\laa\varphi(x,t)\right)dxdt,
\end{equation}
where the velocity $v$ satisfies \eqref{incom} and is given by \eqref{darcy}.
\end{definition}
An analogous definition is considered in the periodic setting
taking $\varphi\in C^{\infty}([0,\tau]\times\TT)$.

\begin{theorem}\label{weak}
Suppose $T_0\in L^2(\RR^N)$ and $0< \alpha\leq 2$. Then, for any $\tau>0$, there exists at least one weak solution
$T\in C([0,\tau];L^2(\RR^N))\cap L^2([0,\tau];H^{\alpha/2}(\RR^N))$ to the DPM equation.
\end{theorem}

\noindent {\bf Proof.} To prove the theorem we modify the system
\eqref{mpv}--\eqref{incom} with a small  viscosity term and we
regularize the initial data.  In particular, for $\ee>0$, we
consider the family $T_\ee$ of solutions given by the system
\begin{equation}\label{w1}
\begin{array}{l}
\displaystyle\frac{\partial T_\ee}{\partial
t}+v_\ee\cdot\nabla T_\ee=-\nu\laa T_\ee+\ee\Delta T_\ee,\\
v_\ee=-\left(\nabla p_\ee+\gamma T_\ee\right),\\
\dv v_\ee=0\\
T_\ee(x,0)=\phi_\ee\ast T_0,
\end{array}
\end{equation}
where $\ast$ denotes the convolution, $\phi_\ee(x)=\ee^{-N}\phi(x/\ee)$ and
$$\phi\in C^\infty_c (\RR^N ), \quad\phi\geq0, \quad\int\phi(x)dx = 1.$$
As we show in the previous section there is a global solution of
\eqref{w1} with $T_\ee\in C([0,\tau];H^s(\RR^N))$ for any $s>0$.
We multiply by $T_\ee$  to get
$$
\frac12\frac d{dt}\left(\|T_\ee\|^2_{L^2}\right) + \nu\|\la^{\alpha/2} T_\ee\|_{L^2}^2\leq 0,
$$
and integrating in time
\begin{equation}\label{w0}
\|T_\ee(\tau)\|^2_{L^2} + 2\nu\int_0^\tau\|\la^{\alpha/2}
T_\ee(s)\|_{L^2}^2ds\leq \|T_0\|^2_{L^2}\qquad \forall \tau.
\end{equation}
In particular we find
\begin{equation}\label{w2}
T_\ee\in C([0,\tau];L^2(\RR^N))\quad\mbox{and}\quad\max\limits_{0\leq t\leq \tau}\|T_\ee(t)\|_{L^2}^2\leq \|T_0\|^2_{L^2}.
\end{equation}

We pass to the limit using the Aubin-Lions compactness lemma (see
\cite{lions}):

\begin{lemma}\label{aubin-lions}
Let $\{f_\ee(t)\}$ be a sequence in $C([0,\tau];H^s(\RR^N))$ such that
\begin{itemize}
\item[$i)$] $\max\{\|f_\ee(t)\|_{H^s}:\: 0\leq t\leq\tau\}\leq C$
\item[$ii)$] for any $\varphi\in C^\infty_c(\RR^{N})$, $\{\varphi f_\ee\}$ is uniformly Lipschitz in the interval of time $[0,\tau]$ with respect to the space $H^r(\RR^N)$ with $r<s$, i.e.,
$$
\|\varphi f_\ee(t_2)-\varphi f_\ee(t_1)\|_{H^r}\leq C_s|t_2-t_1|\qquad 0\leq t_1,t_2\leq\tau.
$$
\end{itemize}
Then, there exists a subsequence  $\{f_{\ee_j}(t)\}$ and $f\in C([0,\tau];H^s(\RR^N))$ such that for all $\lambda\in (r,s)$
$$
\max\limits_{0\leq t\leq \tau}\|\varphi f_{\ee_j}(t)-\varphi f(t)\|_{H^\lambda}\to 0\qquad\hbox{as }j\to\infty.
$$
\end{lemma}

First, by \eqref{w2} we get $T_\ee\in C([0,\tau];L^2(\RR^N))$ and
$i)$ in the space $L^2(\RR^N)$. Next, we prove that the family
$T_\ee$ is Lipschitz in some space $H^{-r}(\RR^N)$ with $r>N/2+2$.
Since  $T_\ee$ is a strong solution of \eqref{w1} and continuous it follows
\begin{equation}\label{w3}
\|\varphi T_\ee(t_2)-\varphi T_\ee(t_1)\|_{H^{-r}}=\left\|\int_{t_1}^{t_2}\varphi \frac d{dt}T_\ee dt\right\|_{H^{-r}}\leq \max\limits_{t_1\leq t\leq t_2}\{A(t)\}|t_2-t_1|,
\end{equation}
where
$$
A(t)=\|\dv(\varphi v_\ee T_\ee)\|_{H^{-r}}+\|\dv(\varphi)v_\ee T_\ee\|_{H^{-r}}+\nu\|\varphi\laa T_\ee\|_{H^{-r}}+\ee\|\varphi\Delta T_\ee\|_{H^{-r}}.
$$
Applying the property that  the Fourier transform of the product
is the convolution of the respective Fourier transforms, we have
$$
\Big|\int_{\RR^N}\hat\varphi(\eta)|\xi-\eta|^\alpha \hat T_\ee(\xi-\eta)d\eta\Big|\leq C \int_{\RR^N}(|\xi|^\alpha+|\eta|^\alpha)|\hat\varphi(\eta)| |\hat T_\ee(\xi-\eta)|d\eta\leq (1+|\xi|^{\alpha})\|\varphi\|_{H^\alpha}\|T_\ee\|_{L^2},
$$
and it yields
$$
\|\varphi\laa T_\ee\|_{H^{-r}}\leq C(\varphi)\|T_\ee\|_{L^2}\Big(\int\limits_{\RR^N}\frac{(1+|\xi|^\alpha)^2 }{(1+|\xi|^2)^r}d\xi\Big)^{1/2}\leq C(r,\varphi)\|T_0\|_{L^2}.
$$
Analogously,
$$
\|\varphi\Delta T_\ee\|_{H^{-r}}\leq C(r,\varphi)\|T_0\|_{L^2}.
$$
We have
$$
\|\!\!\dv(\varphi) v_\ee T_\ee\|_{H^{-r}}\leq C(r)\|\widehat{\dv(\varphi) v_\ee T_\ee}\|_{L^\infty}\!\leq C(r)\|\widehat{\dv(\varphi)}\|_{L^1}\|\widehat{v_\ee T_\ee}\|_{L^\infty}\!\leq C(r,\varphi)\|v_\ee\|_{L^2}\| T_\ee\|_{L^2},
$$
and by \eqref{w2} and the fact that the velocity satisfies \eqref{2.4}, it follows:
$$
\|\!\!\dv(\varphi) v_\ee T_\ee\|_{H^{-r}}\leq C(r,\varphi)\|T_0\|_{L^2}^2.
$$
In a similar way
$$
\|\dv(\varphi v_\ee T_\ee)\|_{H^{-r}}\leq \|\varphi v_\ee T_\ee\|_{H^{1-r}
}\leq C(s,\varphi)\|T_0\|_{L^2}^2.
$$

 From \eqref{w3},  the condition $ii)$ of the Aubin-Lions
lemma is satisfied. Therefore, there exists a subsequence and a function $T\in
C([0,\tau];L^2(\RR^N))$ such that
\begin{equation}\label{w4}
T_\ee\rightharpoonup T \hbox{ in }L^2\hbox{ a.e. } t \hbox{ and }
\max\limits_{0\leq t\leq \tau}\|\varphi T_{\ee}(t)-\varphi T(t)\|_{H^\lambda}\to 0\quad\hbox{as }\lambda\in(-r,0).
\end{equation}

We pass to the limit in the weak formulation of the problem \eqref{w1}, i.e.,
$$
0=\int\limits_{\RR^N}T_\ee(x,0)\varphi(x,0)dx+\int\limits_0^\tau\int\limits_{\RR^N}T_\ee(x,t)
\left(\partial_t\varphi(x,t)+ v_\ee(x,t)\cdot\nabla\varphi(x,t)-\nu\laa\varphi(x,t)+\ee\Delta\varphi(x,t)\right)dxdt,
$$
and we obtain
$$
0=\int\limits_{\RR^N}T_0(x)\varphi(x,0)dx
+\int\limits_0^\tau\int\limits_{\RR^N}T(x,t)
\left(\partial_t\varphi(x,t)-\nu\laa\varphi(x,t)\right)dxdt+\lim_{\ee\to 0}\int\limits_0^\tau\int\limits_{\RR^N}T_\ee (v_\ee\cdot\nabla\varphi) dxdt.
$$
Next, we decompose the non-linear term
$$
\int\limits_0^\tau\int\limits_{\RR^N}T_\ee (v_\ee\cdot\nabla\varphi) dxdt=\int\limits_0^\tau\int\limits_{\RR^N}(T_\ee-T) (v_\ee\cdot\nabla\varphi) dxdt+ \int\limits_0^\tau\int\limits_{\RR^N}T (v_\ee\cdot\nabla\varphi) dxdt.
$$
Using the Fourier transform in the first term we have
\begin{eqnarray*}
\left|\int\limits_0^\tau\int\limits_{\RR^N}(T_\ee-T) (v^\ee\cdot\nabla\varphi) dxdt\right|&\leq & \int\limits_0^\tau\|v_\ee\|_{H^{\alpha/2}}\|(T_\ee-T)\nabla\varphi\|_{H^{-\alpha/2}}dt\\
&\leq & \max\limits_{0\leq t\leq\tau}\|(T_\ee-T)\nabla\varphi\|_{H^{-\alpha/2}}\int\limits_0^\tau\left(\|v_\ee\|_{L^2}+\|\la^{\alpha/2}v_\ee\|_{L^2}\right)dt.
\end{eqnarray*}
Due to \eqref{w0} and \eqref{w2} we get
$$
\int\limits_0^\tau\left(\|v_\ee\|_{L^2}+\|\la^{\alpha/2}v_\ee\|_{L^2}\right)dt\leq c(\tau) \|T_0\|_{L^2}.
$$
Then, by \eqref{w4} we have
$$
\lim_{\ee\to 0}\int\limits_0^\tau\int\limits_{\RR^N}T_\ee (v_\ee\cdot\nabla\varphi) dxdt=\int\limits_0^\tau\int\limits_{\RR^N}T (v\cdot\nabla\varphi) dxdt
$$
and we conclude the proof of theorem~\ref{weak}. \quad\eop

\smallskip

\begin{remark}\label{w-periodic}
Analogous result of theorem~\ref{weak} follows, with a similar argument in the torus $\TT^N$, with periodic boundary conditions.
\end{remark}

We continue this section mentioning the  existence result of weak
solutions obtained for the non-homogeneous equation \eqref{mpv2}.

\begin{theorem}\label{weakf}
Let $\tau > 0$ be arbitrary. Then for every $T_0\in L^2$ and $f\in L^2(0; T;H^{-\frac\alpha 2})$, there exists a weak solution of
(\ref{mpv2}) satisfying $T(x,0)=T_0(x)$ and $T\in C([0,\tau];L^2)\cap L^2([0,\tau];H^{\alpha/2})$.
\end{theorem}

The proof is similar to theorem~\ref{weak}.

\smallskip

Although weak solutions may not be unique, there is at most one solution in the
class of ``strong'' solutions in the sub-critical case. This fact is well known for the quasi-geostrophic equation (see \cite{cons-wu}) and it is a generalization of the classical Leray-Prodi-Serrin condition, related to the uniqueness of the solutions to the 3D Navier-Stokes equation (see \cite{teman}).

\begin{theorem}\label{serrin}
Assume that $\alpha\in(1,2]$, $\tau>0$ and  $T$ a weak solution of DPM with $T_0\in L^2$.
Then, there is an unique weak solution satisfying:
\begin{eqnarray}\label{lpserrin}
T\in C([0,\tau];L^2)\cap L^2([0,\tau];H^{\frac\alpha 2})\cap L^p([0,\tau];L^q),
\end{eqnarray}
for $q>N/(\alpha-1)$, and $p=\alpha/(\alpha-N/q-1).$
\end{theorem}

\smallskip

\noindent {\bf Proof.} We take the difference $T=T^1-T^2$ of two
solutions $T^1$ and $T^2$ of DPM with same initial data.
Considering $v=v^1-v^2$, with $v^1$ and $v^2$ being the velocities
corresponding to $T^1$ and $T^2$, then $T$ satisfies
$$
\frac{\partial T}{\partial
t}+v\cdot\nabla T^1+ v^2\cdot\nabla T +\nu\laa T= 0,
$$
or analogously
$$
\frac{\partial T}{\partial
t}+\dv(v T^1)+\dv(v^2 T)+\nu\laa T= 0.
$$
We multiply  the equation by $\la^{-1}T$ and integrate by parts in
the nonlinear terms to obtain
$$
\frac{d}{dt}\|\la^{-\frac12} T\|_{L^2}^2+ \nu \|\la^\frac\alpha2(\la^{-\frac12}T)\|_{L^2}^2\leq \left|\int_\RR (T\,v_2)\cdot (\nabla(\la^{-1}T))dx\right|
+\left|\int_\RR (T_1v)\cdot(\nabla (\la^{-1}T))dx \right|.
$$
We take $(1/q)+(2/p)=1$ and it yields
$$
\left|\int_\RR (T\,v^2)\cdot (\nabla(\la^{-1}T))dx\right|
\leq \|v^2\|_{L^q}\|T\|_{L^p}\|\nabla(\la^{-1}T)\|_{L^p}
$$
and
$$
\left|\int_\RR (T^1v)\cdot(\nabla (\la^{-1}T))dx\right|\leq \|T^1\|_{L^q}\|v\|_{L^p}\|\nabla(\la^{-1}T)\|_{L^p}.
$$
Since $\nabla (\la^{-1})=(R_1,R_2)$  we obtain
$$
\left|\int_\RR (T\,v^2)\cdot (\nabla(\la^{-1}T))dx\right|+\left|\int_\RR (T^1v)\cdot(\nabla (\la^{-1}T))dx\right|\leq C(\|T^1\|_{L^q}+\|T^2\|_{L^q})\|T\|^2_{L^p}.
$$
The inequality for the Riesz potential \eqref{rieszpotential} gives
$$\|T\|_{L^p}\leq C\|\la^{\frac{N}{2q}}T\|_{L^2}=C\|\la^{\frac{N}{2q}+\frac12}(\la^{-\frac12}T)\|_{L^2}.$$

For $q$ large enough, we can get $(N/2q)+1/2<\alpha/2$, and the
following interpolation inequality
$$
\|\la^sf\|_{L^2}\leq \|f\|^{\gamma}_{L^2}\|\la^{r}f\|^{1-\gamma}_{L^2}
$$ with $s<r$, $0<\gamma<1$, yields
$$
\|\la^{\frac{N}{2q}+\frac12}(\la^{-\frac12}T)\|_{L^2}\leq \|\la^{-\frac12}T\|^{\gamma}_{L^2}\|\la^{\frac\alpha2}(\la^{-\frac12}T)\|^{1-\gamma}_{L^2},
$$
for $\gamma=(\alpha-N/q-1)/\alpha$.
Finally
$$
\frac{d}{dt}\|\la^{-\frac12} T\|_{L^2}^2+ \nu \|\la^\frac\alpha2(\la^{-\frac12}T)\|_{L^2}^2\leq  C(\|T_1\|_{L^q}+\|T_2\|_{L^q})\|\la^{-\frac12}T\|^{2\gamma}_{L^2}
\|\la^{\frac\alpha2}(\la^{-\frac12}T)\|^{2(1-\gamma)}_{L^2}
$$ and
$$
\frac{d}{dt}\|\la^{-\frac12} T\|_{L^2}^2\leq  \frac C\nu (\|T_1\|_{L^q}+\|T_2\|_{L^q})^{1/\gamma}\|\la^{-\frac12}T\|^{2}_{L^2},
$$
which allows to conclude the proof.\quad\eop
\begin{remark}
If we take $T_0\in L^2\cap L^q$ for $q>N/(\alpha-1)$ we can
construct, as before, a solution that satisfies $\|T\|_{L^q}\leq
\|T_0\|_{L^q}$, and in particular we have
$$
T\in L^{\infty}([0,\tau];L^q).
$$
Then this solution is unique in this space.
\end{remark}

\begin{remark}\label{semigroup}
Suppose that $1<\alpha\leq 2$ , $\nu > 0$, $s > 0$, $f\in L^p\cap H^{s-\frac\alpha2}(\TT^N)$ and
$$
T_0 \in H^s\cap L^p(\TT^N) , \quad\hbox{ where } 0\leq \frac 1 p <\frac{\alpha-1}N.$$
Then,  there is a weak solution $T$ of \eqref{mpv2} such that
$$
T\in C([0,\tau];H^s(\TT^N))\cap L^2(0,\tau;H^{s+\frac\alpha 2}(\TT^N)).
$$
The proof follows applying  an analogous analysis as in the proof of theorem~\ref{sub}. Moreover,  as in theorem~\ref{serrin}, the
uniqueness of weak solutions for the non-homogeneous equation \eqref{mpv2} is also obtained.

\end{remark}

\smallskip\

\section{Decay estimates}

Here, we obtain the decay of the solutions of \eqref{mpv}--\eqref{incom}. The key of the argument is the following positivity lemma.

\begin{lemma}\label{positive}
Suposse $\alpha\in[0,2]$, $\Omega=\RR^N, \TT^N$ and $T$, $\la^\alpha T\in L^p(\Omega)$ where $p\geq 2$. Then
$$
\int\limits_\Omega|T|^{p-2}T\la^\alpha Tdx\geq \frac2p\int\limits_\Omega\left(\la^\frac\alpha 2 |T|^\frac p2\right)^2dx.
$$
\end{lemma}

This lemma is a consequence of different versions of the
positivity lemma obtained in \cite{Resnick,cor1,cor2,ju}.

The immediate consequence of the previous lemma is the following decay results
in the $L^p$ space for solutions of (\ref{mpv})--\eqref{incom}:

\begin{corollary}\label{decay}  Suppose that $T_0\in L^p$ where $p\in [2,+\infty)$ and $T$ is a weak solution of (\ref{mpv})--\eqref{incom}.
\begin{enumerate}
\item If $\Omega=\TT^N$ and the mean value of $T_0$ is zero, then
$$
\|T(t)\|_{L^q}\leq \|T_0\|_{L^p}\exp\left(-\frac{2\nu\lambda_1^{\alpha}}p t\right), \qquad q\in[1,p],
$$
where $\lambda_1>0$ is the first positive eigenvalue of $\la $.
\item If $\Omega=\RR^N$, then
$$
\displaystyle\|T(t)\|_{L^p}\leq \|T_0\|_{L^p}\left[1+\displaystyle\frac{4\alpha \nu c \|T_0\|_{L^p}^\frac{2p\alpha}{N(p-2)}t}{N(p-2)\|T_0\|_{L^2}^\frac{2p\alpha}{N(p-2)}} \right]^{\displaystyle-\frac{N(p-2)}{2p\alpha}},
$$
with $c$ depending on $\alpha$ and $N$.
\end{enumerate}
\end{corollary}

\noindent {\bf Proof.} We multiply the equation (\ref{mpv}) by $|T|^{p-2}T$ and applying the lemma~\ref{positive} we obtain
$$
\frac d{dt}\|T(t)\|_{L^p}^p\leq-2\nu\int|\la^\frac{\alpha}{2} T^\frac p2|^2dx.
$$

In the case $\Omega=\TT^N$  we get
$$
\frac d{dt}\|T(t)\|_{L^p}^p+2\nu \lambda_1^\alpha\|T(t)\|_{L^p}^p\leq 0,
$$
and the above inequality gives the exponential decay of $\|T(t)\|_{L^q} $ for $q\in[1,p]$.

In the case $\Omega=\RR^N$, using
Gagliardo-Nirenberg inequality, we have
$$
\frac d{dt}\|T(t)\|_{L^p}^p\leq-2\nu c\left(\int| T|^\frac{pN}{N-\alpha}dx\right)^\frac{N-\alpha}{N},
$$
with $c$ depending on $\alpha$ and $N$. By interpolation  we get
$$
\|T\|_{L^p} \leq \|T\|_{L^2}^{1-\beta} \left(\int| T|^\frac{pN}{N-\alpha}dx\right)^{\beta\frac{N-\alpha}{pN}}
$$
with
$$
\beta=\frac{N(p-2)}{N(p-2)+2\alpha}.
$$
Therefore,
$$
\frac d{dt}\|T(t)\|_{L^p}^p +2c\nu\|T\|_{L^2}^{p-\frac p\beta}\|T\|_{L^p}^\frac p\beta\leq 0,
$$
 since $\beta\in(0,1)$ and $\|T(t)\|_{L^2}\leq \|T_0\|_{L^2}$ yields
$$
\frac d{dt}\|T(t)\|_{L^p}^p +2\nu c\displaystyle\frac{\|T\|_{L^p}^\frac p\beta}{\|T_0\|_{L^2}^{\frac p\beta - p}}\leq 0.
$$
We integrate
$$
\frac{\beta}{\beta-1}\left(\|T(t)\|_{L^p}^{p(1-\frac1\beta)}-\|T_0\|_{L^p}^{p(1-\frac1\beta)}\right)\leq -\displaystyle\frac{2\nu c t}{\|T_0\|_{L^2}^{\frac p\beta - p}}.
$$
Again, since $\beta\in (0,1)$, we have
$$
\left(\displaystyle\frac{\|T(t)\|_{L^p}}{\|T_0\|_{L^p}}\right)^{p(1-\frac1\beta)}\geq 1+ \frac{1-\beta}{\beta}\displaystyle\frac{2\nu c\|T_0\|_{L^p}^{\frac p\beta - p} t}{\|T_0\|_{L^2}^{\frac p\beta - p}},
$$
hence
$$
\|T(t)\|_{L^p}\leq \|T_0\|_{L^p}\left( 1+\frac{1-\beta}\beta\displaystyle\frac{2\nu c \|T_0\|_{L^p}^{\frac p\beta - p} t}{\|T_0\|_{L^2}^{\frac p\beta - p}}\right)^{-\displaystyle\frac\beta{p(1-\beta)}}.
$$
Then, by definition of $\beta$, it follows the polynomial decay. \eop

\begin{remark}\label{l-infty-est}
As a consequence of the previous lemma for the case $\Omega=\RR^N$, we obtain the following estimate for the $L^\infty$-norm:
$$
\|T(t)\|_{L^\infty}\leq \|T_0\|_{L^\infty},
$$
which can be improved  as in \cite{cor2}
$$
\|T(t)\|_{L^\infty}\leq \|T_0\|_{L^\infty}\left(1+ \alpha ct\|T_0\|_{L^\infty}^\alpha\right)^{-\frac{1}{\alpha}}.
$$
\end{remark}

\smallskip\

\section{Long time behavior}

In this section we study some attracting properties of the solutions of \eqref{mpv2} with $\alpha\in (1,2]$.

We introduce an abstract framework for studying the
asymptotic behavior of this system with respect to
two topologies, weak and strong, depending on the uniqueness of the solution. Each such system possesses a global attractor in the
weak topology, but not necessarily in the strong topology, and in general, they are different.

First, we recall some definitions of \cite{teman2}. Let $X$ be a
complete metric space. A {\it semiflow} on $X$,
$\omega:[0,\infty)\times X\to X$, is defined to be a mapping
$\omega(t,x)=S(t)x$ that satisfies the following conditions:
$S(0)x=x$ for all $x\in X$; $\omega$ is continuous; the semigroup
condition, i.e., $S(s)S(t)x=S(s+t)x$ for all $s,t\geq0$ and $x\in
X$, is valid.

A semiflow is {\it point dissipative} if there exists a bounded set $B\subset X$ such that for any $x\in X$ there is a time $\tau(x)$ such that
$\omega(t,x)\in B$ for all $t>\tau(x)$. In this case $B$ is referred as an an {\it absorbing set} for the semiflow $S(t)$.

A semiflow is {\it compact} if for any bounded $B\subset X$ and $t>0$, $S(t)B$ lies in a compact subset in $X$.

$\mathcal{A}\subset X$ is a {\it global attractor} if it satisfies
the following conditions: $\mathcal{A}$ is nonempty, invariant and
compact; $\mathcal{A}$ possesses an open neighborhood
$\mathcal{U}$ such that, for every initial data $u_0$ in
$\mathcal{U}$, $S(t)u_0$ converges to $\mathcal{A}$ as
$t\to\infty$:
$$
{\rm dist}(S(t)u_0,\mathcal{A})\to 0 \quad\hbox{ as } t\to\infty.
$$
Recall that the distance of a point to a set is defined by
$$
{\rm d} (x,\mathcal{A})=\inf\limits_{y\in\mathcal{A}} {\rm d}(x,y).
$$

Now, we state the following result about the theory of global attractors (see \cite{teman2}):

\begin{theorem}\label{teman-atractor}
Let $S(t)$ be a point dissipative, compact semiflow on a complete metric space $X$. Then $S(t)$ has a global attractor in $X$.

Furthermore,
the global attractor attracts all bounded sets in $X$, is the maximal bounded absorbing set and minimal invariant set for the inclusion relation.

Assuming in addition that $X$ is a Banach space, $\mathcal{U}$ is convex and for any $x\in X$, $S(t)x:\RR_+\to X$ is continuous. Then, $\mathcal{A}$ is also connected.
\end{theorem}

In the case that $\mathcal{U}=X$, $\mathcal{A}$ is called the {\it global attractor} of the semigroup $\{S(t)\}_{t\geq 0}$ in $X$.

\smallskip

\subsection{Strong attractor}

 From Remark~\ref{semigroup}, we see immediately that for any
$1<\alpha\leq 2$ and with $s>(N/2)+1-\alpha$, the solution
operator of the porous medium equation (\ref{mpv2}) well defines a
semigroup in the space $H^s$.

We begin this section with some useful a priori estimates of the
solutions \eqref{mpv2} with $f\in L^p$.

\begin{lemma}\label{Lp-esti}
Let $T = T(x,t)$ be a solution of (\ref{mpv2}), the initial data
$T_0\in L^p$ with zero mean value and $p\geq 2$. Then,
$\|T\|_{L^p}$ is uniformly bounded with respect to
$\|T_0\|_{L^p}$. In particular,
\begin{equation}\label{bLp}
\|T(t)\|_{L^p}\leq \left(\|T_0\|_{L^p} -\frac{p}{\nu\lambda_1^\alpha}\|f\|_{L^p} \right)\exp\left\{-\frac{\nu\lambda_1^\alpha}{p}t\right\} +\frac{p}{\nu\lambda_1^\alpha}\|f\|_{L^p},
\end{equation}
and there exists an absorbing ball in $L^p$. Moreover, for $T_0\in L^2$, we get
\begin{equation}\label{bL2}
\nu\int^{t+1}_t\|\la^\frac\alpha2T(s)\|_{L^2}^2ds\leq \left(\|T_0\|_{L^2}^2 -\frac{\|f\|_{L^2}^2}{\nu\lambda_1^\alpha} \right)\exp\left\{-\nu\lambda_1^\alpha t\right\}+\frac{\|f\|_{L^2}^2}{\nu\lambda_1^\alpha} +\frac{1}{\nu}\|\la^{-\frac{\alpha}2} f\|_{L^2}.
\end{equation}
\end{lemma}

\noindent {\bf Proof}. We multiply the equation (\ref{mpv2}) by $p|T|^{p-2}T$ with $p\geq 2$. Integrating, using lemma~\ref{positive} and applying Holder's inequality, we get
\begin{equation}\label{5.3}
\frac d{dt}\|T(t)\|_{L^p}^p+2\nu\|\la^\frac{\alpha}{2} T^\frac p2\|^2_{L^2}\leq p\|f\|_{L^p}\|T\|_{L^p}^{p-1}.
\end{equation}
We denote $\lambda_1$ as the first eigenvalue of $\la$. Since $T$ is mean zero, we have
$$
\|\la^\frac{\alpha}{2} T^\frac p2\|^2_{L^2}\geq \lambda_1^\alpha \|T\|_{L^p}^p.
$$
Therefore,
$$
\frac d{dt}\|T(t)\|_{L^p}+\frac{\nu\lambda_1^\alpha}{p} \|T\|_{L^p}\leq\|f\|_{L^p},
$$
and integrating  we prove the estimate \eqref{bLp}.

Now, to prove \eqref{bL2}, we multiply the equation (\ref{mpv2}) by $2T$ and by \eqref{5.3}
$$
\frac d{dt}\|T(t)\|_{L^2}^2+2\nu\|\la^\frac{\alpha}{2} T\|^2_{L^2}\leq 2\|f\|_{L^2}\|T\|_{L^2}.
$$
Integrating and applying the Young inequality, we get
$$
\frac d{dt}\|T(t)\|_{L^2}^2+\nu\|\la^\frac{\alpha}{2} T\|^2_{L^2}\leq \frac{1}{\nu}\|\la^{-\frac{\alpha}2} f\|_{L^2}^2\leq \frac{1}{\nu}\|f\|_{L^2}^2.
$$
Integrating and using that $\lambda_1$ is the first eigenvalue of $\la$, we obtain
\begin{equation}\label{bL22}
\|T(t)\|_{L^2}^2\leq \left(\|T_0\|_{L^2}^2 -\frac{\|f\|_{L^2}^2}{\nu\lambda_1^\alpha} \right)\exp\left\{-\nu\lambda_1^\alpha t\right\}+\frac{\|f\|_{L^2}^2}{\nu\lambda_1^\alpha}.
\end{equation}
On the other hand, integrating between $t$ and $t+1$, we have
\begin{equation}\label{L2}
\|T(t+1)\|_{L^2}^2+\nu\int^{t+1}_t\|\la^\frac{\alpha}{2} T(s)\|^2_{L^2}ds\leq \|T(t)\|_{L^2}^2+\frac{1}{\nu}\|\la^{-\frac{\alpha}2} f\|_{L^2},
\end{equation}
and we get the estimate \eqref{bL2} by \eqref{bL22}.\quad\eop


\begin{lemma}\label{Hs-esti}
Let $T = T(x,t)$ be a solution of (\ref{mpv2}), $T_0\in H^s$ with
zero mean value and $s>(N/2)+1-\alpha$. Then, $\|\la^sT\|_{L^2}$
is uniformly bounded with respect to $\|\la^sT_0\|_{L^2}$ and
there exists an absorbing ball in the space $H^s$. Moreover, we
have
\begin{equation}\label{blas}
\int^{\tau}_0\|\la^{s+\frac\alpha2}T\|_{L^2}^2dt<+\infty
\end{equation}
and
\begin{equation}\label{blas2}
\nu\int^{t+1}_t\|\la^{s+\frac\alpha2}T\|_{L^2}^2 \qquad\hbox{ is uniformly bounded with respect to $\|\la^sT_0\|_{L^2}$}.
\end{equation}

\end{lemma}

\noindent {\bf Proof}. We have that $\alpha\in ( 1,2)$ and that
$T_0 \in H^s$, where $s>(N/2)+1-\alpha$.  If $s \in
((N/2)+1-\alpha,N)$, let $r = s$. If $s \in [N,+\infty)$, let $r$
be any real number in $((N/2)+1-\alpha,N)$. Then, $T_0\in
H^s\subseteq H^r\subset L^p$, where
$$
\frac 1 p=\frac12-\frac rN<\frac{\alpha-1}N.
$$

We multiply the equation (\ref{mpv2}), with an initial data $T_0$
belonging to $L^p\cap H^a$, by $\la^{2a}T$ with $0<a\leq s$. By an
analogous analysis as in the proof of theorem~\ref{sub} we get
$$
\frac d{dt}\|\la^aT(t)\|_{L^2}^2+\nu\|\la^{a+\frac{\alpha}{2}} T\|^2_{L^2}\leq \frac1\nu\|\la^{a-\frac{\alpha}{2}}f\|_{L^2}
+c \|T\|_{L^{p}}\|\la^{a+\frac12+\frac N{2p}} T\|_{L^2}^2.
$$
Then, by \eqref{bLp}, $T\in L^\infty(0,\infty; L^p)$ and
\begin{eqnarray*}
\frac d{dt}\|\la^aT(t)\|_{L^2}^2+\nu\|\la^{a+\frac{\alpha}{2}} T\|^2_{L^2}&\leq & \frac1\nu\|\la^{a-\frac{\alpha}{2}}f\|_{L^2}^2
+C\|\la^{a+\beta} T\|_{L^2}^2.
\end{eqnarray*}
Now, using the Gagliardo-Nirenberg and Holder inequalities, we have
$$
C\|\la^{a+\beta} T\|_{L^2}^2\leq \frac\nu2 \|\la^{a+\frac{\alpha}{2}} T\|^2_{L^2}+ \frac C\nu\|\la^{a} T\|^2_{L^2}
$$
and
$$
\frac d{dt}\|\la^aT(t)\|_{L^2}^2+\frac\nu2\|\la^{a+\frac{\alpha}{2}} T\|^2_{L^2}\leq  \frac1\nu\|\la^{a-\frac{\alpha}{2}}f\|_{L^2}^2
+\frac C\nu\|\la^{a} T\|^2_{L^2}.
$$
Next, following from the above inequality and (\ref{bL2}) the
uniform boundedness of $\|\la^aT(t)\|_{L^2}$ with respect to
$\|\la^aT_0\|_{L^2}$ can be obtained for $a \leq \alpha/2$ from
applying the Uniform Gronwall lemma (see Remark~\ref{gronwall}).
$$
\|\la^aT(t+1)\|_{L^2}^2\leq \frac{1}{\nu} e^\frac c\nu\left(\left(\|T_0\|_{L^2}^2 -\frac{\|f\|_{L^2}^2}{\nu\lambda_1^\alpha} \right)\exp\left\{-\nu\lambda_1^\alpha t\right\}+\frac{\|f\|_{L^2}^2}{\nu\lambda_1^\alpha} +\frac{1}{\nu}\|\la^{-\frac{\alpha}2} f\|_{L^2}^2+ \|\la^{a-\frac{\alpha}2} f\|_{L^2}^2\right).
$$
This estimate can assure us that also it gives us an absorbing
ball of the solutions in the space $H^a$ with $0<a\leq \alpha/2$.

Moreover, we have \eqref{blas} and \eqref{blas2} for $0<a\leq
\alpha/2$. Therefore, with these estimates and a bootstrapping
argument,  the uniform boundedness of $\|\la^sT(t)\|_{L^2}$ is
indeed valid for any $s>(N/2)+1-\alpha$ by using the Uniform
Gronwall lemma again. This also gives us, as before, an absorbing
set in the space $H^s$ for any $s>(N/2)+1-\alpha$ and the
estimates \eqref{blas} and \eqref{blas2}. \quad\eop

\begin{remark}\label{gronwall}(Uniform Gronwall lemma)
Let $g$, $h$ and $y$ be non-negative locally integrable functions
on $(t_0,\infty)$ such that
$$
\frac{dy}{dt}\leq gy+ h,\qquad \forall t\geq t_0,
$$
 and
$$
\int^{t+r}_t g(s)ds\leq c_1,\quad \int^{t+r}_t h(s)ds\leq
c_2,\quad \int^{t+r}_t y(s)ds\leq c_3,\quad\forall t\geq t_0,
$$
where $r$, $c_1$, $c_2$ and $c_3$ are positive constants. Then,
$$
y(t+r)\leq \left(\frac{c_3}r+c_2\right)e^{a_1}, \quad\forall t\geq
t_0.
$$
The proof of this estimate is shown in \cite{teman2}.
\end{remark}

Now, we prove a condition to apply the
theorem~\ref{teman-atractor}:  the continuity of the solutions of
\eqref{mpv2} in the space $H^s$ with respect to $t$.

\begin{lemma}\label{Hs-conti}
Let $T = T(x,t)$ be a solution of (\ref{mpv2}), $T_0\in H^s$ with
zero mean value and $s>(N/2)+1-\alpha$. Then, $\la^s T \in
C(0,\tau; L^2)$.
\end{lemma}

\noindent {\bf Proof}. By lemma~\ref{Hs-esti} we have that $\la^s
T \in L^2(0,\tau; H^\frac\alpha2)$. According to the Aubin-Lions
compactness results (see \cite{simon}) we just need to show that
$\la^s T_t \in L^2(0,\tau; H^{-\frac\alpha2})$. From the equation
(\ref{mpv2}) we get
\begin{equation}\label{bL2t}
\|\la^sT_t\|_{H^{-\frac\alpha2}}\leq \|\la^{1+s-\frac\alpha2}(T v)\|_{L^2} +\|\la^{s+\frac\alpha2}T\|_{L^2}+\|\la^{s-\frac\alpha2} f\|_{L^2}.
\end{equation}
Using \eqref{resnisk} and the integral formulation of the velocity
we have
$$
\|\la^{1+s-\frac\alpha2}(T v)\|_{L^2}\leq C\|T\|_{L^p} \|\la^{1+s-\frac\alpha2}T\|_{L^q}
$$
where (1/p)+(1/q)=1/2. Now, as in the proof of
lemma~\ref{Hs-esti}, we take $r\leq s$ such that
 $T_0\in H^s\subseteq H^r\subset L^p$, with
$$
\frac 1 p=\frac12-\frac rN<\frac{\alpha-1}N.
$$
Then, considering $q=N/r$, $T\in L^\infty(0,\infty; L^p)$. Since
$r> (N/2)+1-\alpha$ (see the proof of lemma~\ref{Hs-esti}) we have
that
$$
q=\frac{N}{r}<\frac{2N}{N+2-2\alpha}=q^*.
$$
Therefore,
$$
\|\la^{1+s-\frac\alpha2}T\|_{L^q}\leq c\|\la^{1+s-\frac\alpha2}T\|_{L^{q^*}}\leq C \|\la^{s+\frac\alpha2}T\|_{L^2},
$$
applying the Gagliardo-Nirenberg inequality. Thus, coming to \eqref{bL22}, we get
$$
\|\la^sT_t\|_{H^{-\frac\alpha2}}\leq  (C\|T\|_{L^p}+1)\|\la^{s+\frac\alpha2}T\|_{L^2}+\|\la^{s-\frac\alpha2} f\|_{L^2}.
$$
Finally, by \eqref{bLp} and \eqref{blas}, we obtain
$$
\int_0^\tau\|\la^sT_t(t)\|_{H^{-\frac\alpha2}}dt<\infty,
$$
and we conclude the proof.\quad\eop

\begin{lemma}\label{t-conti}
Assume that the initial data, of a solution of equation
(\ref{mpv2}), belongs to $H^s$ with zero mean value and $s>
N/2+1-\alpha$. Then, for any fixed $t > 0$, the solution operator
$S(t)$ is a continuous map from $H^s$ into itself.
\end{lemma}

\noindent {\bf Proof}. We consider  two solutions $T^{(1)}$ and $T^{(2)}$ of the porous medium equation
 (\ref{mpv2}) with two initial data $T^{(1)}_0$ and $T^{(2)}_0$ and velocities $v^{(1)}$ and $v^{(2)}$, respectively.
  Let $T=T^{(1)}-T^{(2)}$ and $v=v^{(1)}-v^{(2)}$. Then, since $\dv(v)=0$, we have for any $\varphi\in H^\frac{\alpha}2$ that
\begin{equation}\label{1L3}
(T_t,\varphi) +\nu (\la^\frac\alpha2 T,\la^\frac\alpha2\varphi)=-(v\cdot\nabla T^{(2)},\varphi)-(v^{(1)}\nabla T,\varphi).
\end{equation}

Setting $\varphi=T$, using that $\dv(v^{(1)})=0$ and Young inequality, we get
$$
\frac12\frac d{dt}\|T\|^2_{L^2}+\nu\|\la^\frac\alpha2T\|_{L^2}^2\leq \|\la T^{(2)}\|_{L^{q_1}}\|T\|_{L^{q_2}}^2
$$
such that $(1/q_1)+(2/q_2)=1$. By Gagliardo-Nirenberg inequality, we obtain
$$
\frac12\frac d{dt}\|T\|^2_{L^2}+\nu\|\la^\frac\alpha2T\|_{L^2}^2\leq \|\la T^{(2)}\|_{L^{q_1}}\|T\|_{L^{2}}^{2(1-a)}\|\la^\frac\alpha 2 T\|_{L^{2}}^{2a}
$$
with $a=N/(q_1\alpha)$, where we will be choosing $q_i$ such that
$a\in(0,1)$. We use again  Young inequality and we have
$$
\frac12\frac d{dt}\|T\|^2_{L^2}+\frac{\nu}{2}\|\la^\frac\alpha2T\|_{L^2}^2\leq c(\nu)\|\la T^{(2)}\|_{L^{q_1}}^{q}\|T\|_{L^{2}}^{2},
$$
denoting $q=1/(1-a)$. Thus, by the Gronwall inequality, it follows
$$
\|T(t)\|^2_{L^2}\leq C(\nu)\|T_0^{(1)}-T_0^{(2)}\|_{L^{2}}^{2}\exp\left(\int_0^t\|\la T^{(2)}(s)\|_{L^{q_1}}^{q}ds\right).
$$
If $s\in ((N/2)+1-\alpha, (N/2)+1-(\alpha/2))$, then we take $r =
s$.  If $s \in [(N/2)+1-(\alpha/2),+\infty)$, we take $r$ any
number in $((N/2)+1-\alpha, (N/2)+1-(\alpha/2))$. Then $H^s
\subseteq H^r$. We choose
$$
q_1=\frac{2N}{2+N-2r-\alpha}>1,
$$
then
$$
a=\frac{2+N-2r-\alpha}{2\alpha}\in (0,\frac 12)\quad \hbox{ and }\quad q<2.
$$
Therefore, using the following Sobolev inclusions
$$
L^q(0,\tau; W^{1,q_1})\subset L^q(0,\tau; H^{r+\frac\alpha2})\subset L^2(0,\tau; H^{r+\frac\alpha2})\subset L^2(0,\tau; H^{s+\frac\alpha2}),
$$
we conclude that
$$
\int^\tau_0\|\la^\frac\alpha2T(t)\|_{L^2}^2dt \leq C(T^{(2)},\tau)\|T_0^{(1)}-T_0^{(2)}\|_{L^{2}}^{2}
$$
Thus, using the Riesz lemma, it is immediate that the solution
operator $S(t)$ is a continuous map from $H^s$ into itself  when $s\in ((N/2)+1-\alpha, \alpha/2]$.

We finish the proof studying the case $s > \alpha/2$. We do so by
checking directly the Lipschitz continuity of the solution
operator in the space $H^s$. We consider $\varphi=\la^{2s}T$ in
\eqref{1L3}, then
\begin{equation}\label{1L4}
\frac12\frac d{dt}\|\la^sT\|^2_{L^2}+\nu\|\la^{s+\frac\alpha2}T\|_{L^2}^2=-(v\cdot\nabla T^{(2)},\la^{2s}T)-(v^{(1)}\nabla T,\la^{2s}T).
\end{equation}
We estimate the two terms on the right-hand side of the
variational formula separately. First, we get
$$
|(v\nabla T^{(2)},\la^{2s}T)|\leq
c(\nu)\|\la^{s\frac{\alpha}2}(v\cdot\nabla
T^{(2)})\|_{L^2}^2+\frac\nu 4\| \la^{s+\frac{\alpha}2}T\|_{L^2}^2.
$$
Using similar estimate as \eqref{resnisk} of Kenig, Ponce and Vega (see \cite{kpv}), we have
$$
\|\la^{s-\frac{\alpha}2}(v^{(1)}\nabla T)\|_{L^2}\leq \|\la^{s-\frac{\alpha}2}T\|_{L^{p_1}}\|\la T^{(2)}\|_{L^{p_2}}+\|T\|_{L^{q_1}}\|\la^{s+1-\frac{\alpha}2} T^{(2)}\|_{L^{q_2}}
$$
with $(1/p_1)+(1/p_2)=1/2$ and $(1/q_1)+(1/q_2)=1/2$. We select
$$
p_1=\frac{2N}{N-\alpha},\qquad p_2=\frac{2N}\alpha,\qquad q_1=\frac{N}{\alpha-1},\qquad q_2=\frac{2N}{N+2-2\alpha}
$$
and, using the Sobolev inequalities yields the following estimate
\begin{equation}\label{1L5}
|(v^{(1)}\nabla T,\la^{2s}T)|\leq c(\nu)\|\la^{s}T\|_{L^2}^2\|\la^{s+\frac{\alpha}2} T^{(2)}\|_{L^2}^2+\frac\nu 4\| \la^{s+\frac{\alpha}2}T\|_{L^2}^2.
\end{equation}
We estimate the other term on the right side of \eqref{1L4}. Since
$(v^{(1)}\cdot\nabla \la^s T, \la^s T)=0$, $\la^s$ and $\nabla$
are commutable we have
$$
|(v^{(1)}\nabla T,\la^{2s}T)|\leq \|\la^s(v^{(1)}\cdot\nabla T)- v^{(1)}\cdot(\la^s\nabla T)\|_{L^2}\|\la^{s}T\|_{L^2}.
$$
Using  estimate of Kenig, Ponce and Vega (see \cite{kpv}) we have
$$
\|\la^s(v^{(1)}\cdot\nabla T)- v^{(1)}\cdot(\la^s\nabla T)\|_{L^2}\leq \|\la v^{(1)}\|_{L^{p_1}} \|\la^{s} T\|_{L^{p_2}} +\|\la^s v^{(1)}\|_{L^{q_1}} \|\la T\|_{L^{q_2}},
$$
with $(1/p_1)+(1/p_2)=1/2$ and $(1/q_1)+(1/q_2)=1/2$. We take
$$
p_1=q_2=\frac{2N}\alpha,\qquad p_2=q_1\frac{2N}{N-\alpha},
$$
and using the Sobolev inequalities we get
\begin{eqnarray*}
|(v^{(1)}\nabla T,\la^{2s}T)|&\leq& \|\la^{s+\frac\alpha2} T^{(1)}\|_{L^{2}} \|\la^{s+\frac\alpha2} T\|_{L^{2}}\|\la^{s}T\|_{L^2}\\
&\leq& c(\nu)\|\la^{s+\frac\alpha2} T^{(1)}\|_{L^{2}}^2 \|\la^{s}T\|_{L^2}^2+\frac{\nu}4 \|\la^{s+\frac\alpha2} T\|_{L^{2}}^2.
\end{eqnarray*}
Therefore, considering this estimate and \eqref{1L5} in \eqref{1L4}, we obtain
$$
\frac d{dt}\|\la^sT\|^2_{L^2}+\nu\|\la^{s+\frac\alpha2}T\|_{L^2}^2\leq
 c(\nu)\left(\|\la^{s+\frac\alpha2} T^{(1)}\|_{L^{2}}^2 + \|\la^{s+\frac\alpha2} T^{(2)}\|_{L^{2}}^2\right)\|\la^{s}T\|_{L^2}^2.
$$
So, by Gronwall's lemma and since
$$
\int^\tau_0 \left(\|\la^{s+\frac\alpha2} T^{(1)}(t)\|_{L^{2}}^2 + \|\la^{s+\frac\alpha2} T^{(2)}(t)\|_{L^{2}}^2\right)dt <\infty,
$$
we get
$$
\|\la^sT(t)\|_{L^2}\leq C(\nu, \|\la^{s+\frac\alpha2} T^{(1)}_0\|_{L^{2}},\|\la^{s+\frac\alpha2} T^{(1)}_0\|_{L^{2}}) \|\la^{s} (T^{(1)}_0-T^{(2)}_0)\|_{L^{2}}
$$
and conclude the proof of the lemma. \quad\eop

\smallskip

Finally, we present the
existence of the global classic  attractor.

\begin{theorem}\label{s-atractor} Let $\alpha\in(1,2]$, $\nu>0$, $s>(N/2)+1-\alpha$ and $f\in H^{s-\alpha}\cap L^p$ time-independent external source.
Then, the operator $S$, such that $S(t)T_0=T(t)$ for any $t>0$ and $T$ solution of \eqref{mpv2}, defines a semigroup in the space $H^s$ and satisfies:
\begin{itemize}
\item [i)] For any $t>0$, $S(t)$ is a continuous compact operator in $H^s$.

\item [ii)] For any $T_0\in H^s$, $S$ is a continuous map from $[0,t]$ into $H^s$.

\item [iii)] $\{S(t)\}_{t\geq0}$ possesses an attractor ${\mathcal A}$ that is compact, connected and maximal in $H^s$. ${\mathcal A}$ attracts all bounded subsets of $H^r$  in the norm of $H^r$, for any $r>\alpha-(N/2)-1$.

\item [iv)] If $\alpha>(N+2)/4$, ${\mathcal A}$ attracts all bounded subsets of periodic functions of $L^2$ in the norm of $H^r$, for any $r>\alpha-(N/2)-1$.
\end{itemize}
\end{theorem}

\noindent {\bf Proof}. Items i) and ii) are already proven in
lemma~\ref{t-conti} and lemma~\ref{Hs-conti}, respectively.

To verify the rest of  items we use results of semigroups and the
existence of their attractors (see theorem~\ref{teman-atractor}).
In particular, we need to prove the existence of a bounded subset
$B_0\subset H^s$, an open subset $U$ of $H^s$, such that
$B_0\subseteq U\subseteq H^s$, and $B_0$ is absorbing in $U$, i.e.
for any bounded subset $B\subset U$, there is a $\tau(B)$, such
that $S(t)B\subset B_0$ for all $t>\tau(B)$. This fact is proved
in lemma~\ref{Hs-conti}.

Item iv) can be checked easily since for $\alpha > (N+2)/4$ we
have that $\alpha > (N/2)+1- \alpha$ and that for any $T_0\in
L^2$,
$$
\int^\tau_0\|\la^\frac\alpha2 T\|_{L^2}< \infty\qquad \forall
\tau>0.
$$
Which together with \eqref{L2} we can conclude the proof.\quad
\eop

\smallskip\

\subsection{Weak atractor}\label{s5.1}

We obtain the following estimate for the time derivative of a
solution of \eqref{mpv2}.

\begin{proposition}\label{deri-esti}
Let $T$ be the weak solution of \eqref{mpv2} obtained in theorem~\ref{weakf}. Then,
$$
\frac{\partial T}{\partial t}\in L^r_{loc}(0,\infty;H^{-\sigma})\:\hbox{ with }\frac N2 +1=\sigma+\frac\alpha r,\quad\sigma\in(0,1).
$$
\end{proposition}

\noindent {\bf Proof.} For a smooth $\varphi$ and $0<a<1$ we have
$$
\left|\,\int\limits_{\TT^N}\dv(v\cdot T) \varphi\right| \leq \sum_{i}\int\limits_{\TT^N}|\la^a(R_j v^j T)| |\la^{1-a}\varphi|,
$$
where $R_j$ are the Riesz transforms. Using the inequality \eqref{resnisk} with (1/p)+(1/q)=(1/2) and $v$ satisfies that
$\|\la^a v\|_{L^r}\leq c\|\la^a T\|_{L^r}$, for  $1< r<\infty$, we obtain
\begin{eqnarray*}
\left|\,\int\limits_{\RR^N}\dv(v\cdot T) \varphi\right| &\leq& c \|\la^a T \|_{L^p}\|T\|_{L^q} \|\la^{1-a} \varphi \|_{L^2}.
\end{eqnarray*}
 By interpolation we get
$$
\|T\|_{L^q}\leq \|T\|_{L^2}^{1-s}\|T\|_{L^\frac{2N}{N-\alpha}}^s\: \hbox{ such that } s=\frac N\alpha\left(1-\frac2q\right).
$$
Since $H^{\alpha/2}\subset W^{a,p}$ and
$$
p= \frac{2N}{N+2a-\alpha},
$$
we have
$$
\|\la^a T\|_{L^p}\leq c\|\la^\frac\alpha2 T\|_{L^2} \:\hbox{ and }\|T\|_{L^\frac{2N}{N-\alpha}}\leq c \|\la^\frac\alpha2 T\|_{L^2}.
$$
Therefore, we obtain
$$
\left|\int\limits_{\RR^N}\dv(v\cdot T) \varphi\right|\leq c \|T\|_{L^2}^{1-s}\|\la^\frac\alpha2 T\|_{L^2}^{1+s} \|\la^{1-a} \varphi \|_{L^2}.
$$
Furthermore, using that weak solutions belong to
$L^\infty(0,\tau;L^2)\cap L^2(0,\tau;H^\frac\alpha2)$ by
theorem~\ref{weakf}, it follows
$$
\frac{\partial T}{\partial t}\in L^\frac2{1+s}_{loc}(0,\infty,H^{a-1}).
$$
We conclude the proof defining $\sigma=2/(1+s)$, $\sigma=1-a$ and
using the relations of $p$, $q$, $a$ and $s$.\quad\eop

\smallskip

By theorem~\ref{weakf} we have that $T$ is a weak solution of
\eqref{mpv2} that satisfies
\begin{eqnarray}\label{gw1}
&&T\in L^2_{loc}[0,\infty; L^2)\cap  L^\infty_{loc}(0,\infty; L^2)\cap L^2_{loc}(0,\infty; H^\frac\alpha2),\\
\label{gw2}&&
\frac{\partial T}{\partial t}\in L^r_{loc}[0,\infty;H^{-\sigma}),\quad\hbox{ with } \frac N2 +1=\sigma+\frac\alpha r,\quad\sigma\in(0,1).
\end{eqnarray}

We observe that $L^2_{loc}[0,\infty; L^2)$ is a complete space
metric. However, the weak solutions is not closed in this space.
We define the space of generalized weak solutions $GW(f)$ formed
by the functions  $T\in L^2_{loc}[0,\infty; L^2)$ such that are
{\it generalized weak} solutions of \eqref{mpv2} ($T\in GW(f)$) if
$T$ satisfies \eqref{mpv2} in the sense of distributions and
satisfies \eqref{gw1} and \eqref{gw2}.

Given $T_1, T_2\in L^p(0,\infty;X)$ ($X$ is a Banach space), we consider the metric
\begin{equation}\label{met}
{\rm d}(T_1,T_2)=\sum\limits_{n=0}^\infty\frac1{2^n}\min\left(1,\|T_1-T_2\|_{L^p(n,n+1;X)}\right)
\end{equation}
for the set $GW(f)$ (similarly for $p=\infty$). This metric is
invariant on $L^2_{loc}[0,\infty; L^2)$ (see \cite{ber}). Then
applying classical compactness results (see \cite{simon}) for each
$f\in L^\infty(0,\infty, L^2)$, the set $GW(f)$ is a closed subset
of $L^2_{loc}[0,\infty,L^2)$. Now, we can state the following
lemma:

\begin{lemma} Given $f\in L^\infty(0,\infty, L^2)$ and consider the space $GW(f)$ with the metric defined \eqref{met}, we have
\item[i)] The set $\{(T,f)\}$ with $T\in GW(f)$ and $\|f\|\leq K_0$ in the norm of the space $ L^\infty(0,\infty; L^2)$, is closed in
$L^2_{loc}[0,\infty; L^2)\times L^\infty(0,\infty; L^2)$.

\item[ii)] The mapping $S(t)$ is a semiflow on $L^2_{loc}[0,\infty; L^2)$ and $GW(f)$ is a positively invariant subset.

\item[iii)] $S(t)$ restricted to  $GW(f)$ is compact for $t>0$.

\item[iv)] $S(t)$ restricted to  $GW(f)$ is point dissipative.
\end{lemma}

\noindent {\bf Proof}.

 We verify i)  considering a sequence $\{T_k,f_k\}\in L^2_{loc}[0,\infty; L^2)\times L^\infty(0,\infty; L^2)$
  such that $T_n\in GW(f_n)$ and $\|f_n\|\leq K_0$. We have that $f_n\to f$ in $L^\infty(0,\infty;L^2)$
  and $T_n\to T$ in $L^2_{loc}[0,\infty; L^2)$. From the estimates of the time derivative  of $T_n$ and
   by the classical compactness results of \cite{simon}, we easily obtain that $T\in GW(f)$ by using the weak formulation of the solutions.
   Moreover, $\{(T,f)\}$ is closed.

Next, we use a smoothing argument and the bound on the time
derivative of $T$ as in \cite{sell}  to show the continuity of the
semiflow. The positively invariant of $GW(f)$ is immediate.

The items iii) and iv)  are a consequence of the a-priori estimates  and
the compactness result of \cite{simon}.\quad\eop

\smallskip

By the above lemma and using the existence results of global attractors for a point dissipative compact semiflows on a complete metric space (see theorem~\ref{teman-atractor}), we prove the following attractor result:

\begin{theorem}\label{w-atractor} Let $f\in L^2$ independent of $t$ and $\alpha\in(0,2)$. Then, there exists a global attractor
$\mathcal{A}$, subset of the weak solutions of \eqref{mpv2}, for the semiflow generated by the time-shift
on the space of generalized weak solutions $GW(f)$. Moreover, $\mathcal{A}$ attracts all bounded sets in $GW(f)$.
\end{theorem}

We note that the weak attractor
 is defined in a very weak sense and it gives us less useful information than the global attractor in the classic sense.

\begin{remark}\label{time-dep}
In the case of time dependent external source $f\in L^2_{loc}(0,\infty; L^2)$ it is possible to extend the results of the previous theorem.

We define $f_t(\tau)=f(t+\tau)$ is the time-shift of $f$ and we consider the hull $\mathcal{H}^+(\mathcal{F})$ of the positive time translates $f_t$ with $t\geq 0$ of the external force $f\in \mathcal{F}$ where $\mathcal{F}\subset L^2_{loc}(0,\infty; L^2)$ is a set bounded.
Assuming that $\mathcal{H}^+(\mathcal{F})$ is compact with respect to the weak topology of  $ L^2_{loc}(0,\infty; L^2)$, then there exists a global weak attractor $\mathcal{A}$, subset of the weak solutions of \eqref{mpv2}, for the semiflow $S(t)(T,f)=(T_t,f_t)$.
For more details  one can refer to  the works \cite{ber,sell}.
\end{remark}

\smallskip\


\section{Solutions with infinite energy} \label{s6}


For a divergence free velocity field  there exists a
stream function $\psi$, in the two dimensional case, such that
$$v=\nabla^\perp \psi=(-\pa_{x_2}\psi,\pa_{x_1}\psi).$$
We shall choose a stream function of the form
$$\psi(x_1,x_2,t)=x_2f(x_1,t).$$
Taking the rotational over the equation (\ref{darcy}) we obtain
$$
\nabla\times v=\pa_{x_2}v_1-\pa_{x_1}v_2=-\Delta
\psi=-x_2\pa^2_{x_1}f=-\pa_{x_1}T.
$$
Therefore, the function $T$ has the following expression
$$
T(x_1,x_2,t)=x_2\pa_{x_1}f(x_1,t)+\hat{g}(x_2,t)
$$
where we choose
$$
\hat{g}(x_2,t)=\frac{1}{\pi}x_2\int_0^t||\pa_{x_1}f(\tau)||_{L^2(-\pi,\pi)}^2d\tau.
$$

Substituting the expression in (\ref{mpv}), without diffusion in
the two dimensional case, we obtain
\begin{equation}
\pa_tf_{x}=-\pa_tg-ff_{xx}+(f_x)^2+gf_x\label{ecuacion},
\end{equation}
(here and in the sequel of the section, we denote with subscript the derivatives with respect to $x$) where $g$ satisfies
\begin{equation}\label{g}
g(t)=\frac{1}{\pi}\int_0^t||f_x(\tau)||_{L^2(-\pi,\pi)}^2d\tau
\end{equation}
and we define $f$ as
\begin{equation}\label{f}
f(x,t)=\int_{-\pi}^xf_x(x',t)dx'.
\end{equation}
Notice that the
difference between this system  and the one obtained in
\cite{cgo} is that this  has the property of conserving the
mean zero value of the initial data. Indeed, integrating equation
(\ref{ecuacion}) over the interval $[-\pi,\pi)$ and imposing
periodical condition on $f_x$, we have
$$
\pa_t \int_{-\pi}^\pi f_x(x',t)dx'=\left(-f_x(\pi,t)+g(t)\right)\int_{-\pi}^\pi
f_x(x')dx'.
$$
Therefore, if $f_x$ is a solution of equation (\ref{ecuacion}) and
$$\int_{-\pi}^\pi f_x(x',0)dx'=0,$$ we obtain
$$\int_{-\pi}^\pi f_x(x',t)dx'=0, \qquad \forall t>0.$$


\smallskip\

\subsection{Existence.}

In this section we  prove the following theorem.

\begin{theorem}\label{teorema1}
Let $\varphi_0\in H^2(\ff{T})$ with mean zero value and $$M_0=\max
_{x\in \ff{T}}\varphi_0(x).$$ Then, there exist a solution $f(x,t)$ of
the equation (\ref{ecuacion}) with initial datum $f_{x}(x,0)=\varphi_0(x)$ such that $$f_{x}\in
C([0,T),H^2(\ff{T})),$$ with $T=M_0^{-1}$.

\end{theorem}

In order to prove this theorem, first, we add other diffusion term  to the equation \eqref{ecuacion}. Thus, we have the following system
\begin{equation}\label{aux}
\left\{\begin{array}{l}
\pa_tf_{x}=-\pa_tg-ff_{xx}+(f_x)^2+gf_x+\nu(||f_{xx}||_{L^2}^2+g^2)f_{xxx},\\
f_{x}(x,0)=\varphi_0(x),
\end{array}\right.
\end{equation}
where $g$ satisfies \eqref{g} and $\nu>0$. In the next lemma, we prove the global
existence of the solutions of \eqref{aux}.

\begin{lemma}\label{lemaaux}
Let $\varphi_0\in H^3(\ff{T})$ with mean zero value and $\nu>0$.
Then, there exist a function $f(x,t)$ defined by \eqref{f} where $f_x$ is solution of the equation (\ref{aux}) such that $$f_x\in C([0,\infty),H^3(\ff{T})).$$
\end{lemma}

\noindent {\bf Proof.}  We note that if $\varphi_0$ has mean zero value then
$f$ has mean zero value. Multiplying the equation (\ref{aux}) by $f_x$
and integrating over the interval $[-\pi,\pi)$, we obtain
$$
\frac{1}{2}\frac{d}{dt}||f_x||_{L^2}^2=\frac{3}{2}\int_{-\pi}^\pi(f_x)^3+g||f_x||_{L^2}^2-\nu(||f_{xx}||_{L^2}^2+g^2)||f_{xx}||_{L^2}^2.
$$
Therefore,
$$
\frac{1}{2}\frac{d}{dt}||f_x||_{L^2}^2\leq
\frac{3}{2}||f_x||_{L^\infty}||f_{x}||_{L^2}^2+g||f_x||_{L^2}^2-\nu(||f_{xx}||_{L^2}^2+g^2)||f_{xx}||_{L^2}^2.
$$
Using Gagliardo-Niremberg  and Poincar\'e inequalities,  we
have
$$||f_x||_{L^\infty}\leq C||f_x||_{L^2}^\frac{1}{2}||f_{xx}||_{L^2}^\frac{1}{2}\leq C||f_{xx}||_{L^2}.$$
Hence,
$$
\frac{1}{2}\frac{d}{dt}||f_x||_{L^2}^2\leq
C(||f_{xx}||_{L^2}^2||f_{x}||_{L^2}+g||f_x||_{L^2}||f_{xx}||_{L^2})
$$
$$-\nu||f_{xx}||_{L^2}^4-\nu g^2||f_{xx}||_{L^2}^2.$$
And using the Young's inequality we obtain
$$
\frac{1}{2}\frac{d}{dt}||f_x||_{L^2}^2+\frac{\nu}{4}(||f_{xx}||_{L^2}^4+g^2||f_{xx}||_{L^2}^2)\leq
C_{\nu}||f_{x}||_{L^2}^2.
$$
Therefore,
\begin{equation}
||f_x||_{L^2}\leq ||\varphi_0||_{L^2}\exp(C_\nu t),
\end{equation}
and
\begin{equation}
\int_0^T||f_{xx}||_{L^2}^4dt \leq C(||\varphi_0||_{L^2},\nu, T), \qquad \forall T>0.
\end{equation}
Taking a derivative over equation
(\ref{aux}), multiplying by $f_{xx}$ and integrating over the
interval $[-\pi,\pi)$ yield
$$
\frac{1}{2}\frac{d}{dt}||f_{xx}||_{L^2}^2+\nu(||f_{xx}||_{L^2}^2+g^2)||f_{xxx}||_{L^2}^2=\frac{3}{2}\int_{-\pi}^{\pi}f_x(f_{xx})^2dx+g||f_{xx}||_{L^2}^2.
$$
Hence,
$$
\frac{1}{2}\frac{d}{dt}||f_{xx}||_{L^2}^2\leq
\frac{3}{2}||f_{x}||_{L^\infty}||f_{xx}||_{L^2}^2+g||f_{xx}||_{L^2}^2\leq
C(||f_{xx}||_{L^2}^3+g||f_{xx}||_{L^2}^2).
$$
Integrating between  0 and $T$ we obtain
\begin{eqnarray*}
||f_{xx}||_{L^2}^2&\leq&
||\varphi_{0,x}||_{L^2}^2+C(\int_0^T||f_{xx}||_{L^2}^3dt+\int_0^Tg||f_{xx}||_{L^2}^2dt)\\
&\leq& ||\varphi_{0,xx}||_{L^2}^2+ C(T)\int_0^T||f_{xx}||_{L^2}^4,
\end{eqnarray*}
and we can conclude that $||f_{xx}||_{L^2}$ is bounded for all $T<\infty$.

Finally, we estimate $||f_{xxx}||_{L^2}$ and
$||f_{xxxx}||_{L^2}$. Taking two derivatives on the equation
(\ref{aux}), multiplying by $f_{xxx}$ and integrating over the
interval $[-\pi,\pi]$ yield
\begin{eqnarray*}
\frac{1}{2}\frac{d}{dt}||f_{xxx}||_{L^2}^2+\nu(||f_{xx}||_{L^2}^2+g^2)||f_{xxxx}||_{L^2}^2 &=&
\frac{1}{2}\int_{-\pi}^\pi f_x(f_{xxx})^2+g||f_{xxx}||_{L^2}^2\\
&\leq& C(T)||f_{xxx}||_{L^2}^2.
\end{eqnarray*}
Applying Gronwall inequality we have that $||f_{xxx}||_{L^2}$ is
bounded for all $T<\infty$. We obtain that $||f_{xxxx}||_{L^2}$ is
bounded in a similar form and we conclude the proof.\quad\eop

\smallskip

In order to prove theorem~\ref{teorema1} we show some estimates,
independent of $\nu$, of the global solutions of the equation
(\ref{aux}) by the lemma~\ref{lemaaux}, which allows us to
obtain the local existence for the equation (\ref{ecuacion}). Next, we
 prove the following lemma.

\begin{lemma}\label{radamager}
Let $f_x$ be a global  solution of the equation (\ref{aux}) with initial data $\varphi_0$ and $M(t)$ the maximum of $f_x$. Then
\begin{equation}
M(t)+g(t)\leq\frac{M(0)}{1-M(0)t}.
\end{equation}
\end{lemma}

\noindent {\bf Proof.} In this proof we use the techniques of article
\cite{cgo}  for the control of the maximum of the solution of
equation (\ref{ecuacion}). Let denote $x_M(t)$ to be the point
where $f_{x}$ reaches the maximum, then
\begin{eqnarray*}
(M+g)_t&=&M^2+gM+\nu(||f_{xx}||^2+g^2)f_{xx}(x_M(t),t)\\
&\leq& M^2+gM\leq (M+g)^2.
\end{eqnarray*}
Since $g(0)=0$ we obtain
$$
M+g\leq \frac{M(0)}{1-M(0)t},
$$
and the proof is finished.\quad\eop

\smallskip\

\noindent {\bf Proof of  theorem~\ref{teorema1}.} Multiplying the equation (\ref{aux}) by $f_x$ and integrating over the
interval $[-\pi,\pi)$ yields
$$
\frac{1}{2}\frac{d}{dt}||f_{x}||^2_ {L^2}\leq
(M+g)||f_x||^2_{L^2}.
$$
Applying lemma~\ref{radamager} and Gronwall inequality we
have that $||f_x||_{L^2}$ is bounded for all $T<M(0)^{-1}$. In a
similar way, we can obtain that  $||f_{xx}||_{L^2}$ and
$||f_{xxx}||_{L^2}$ are bounded for all $T<M(0)^{-1}$
independently of $\nu$ .

 To finish the proof, we  consider a sequence of solutions
$\{f^\epsilon\}_{\epsilon>0}$ of the equations
\begin{equation}\label{aux2}
\left\{\begin{array}{l}
\pa_t f^\ep_x=-g^\ep_t-f^\ep f_{xx}^\ep+(f^\ep_x)^2+g^\ep f_x^\ep+\ep(||f_{xx}^\ep||_{L^2}+(g^\ep)^2)f_{xx}^\ep\\
f_x^\ep(x,0)=\varphi_0^\ep(x),
\end{array}\right.
\end{equation}
where $\{\varphi_0^\ep\}_{\ep}$ is a sequence in $H^3(\ff{T})$ such
that
$$\lim_{\ep\rightarrow 0}\varphi_0^\ep=\varphi_0\in H^2,$$

$$M^\ep(0)\equiv \max_{x\in \ff{T}}f^\ep_{x0}(x)\leq M(0)\equiv \max_{x\in \ff{T}}\varphi_0(x),$$
and
$$||\varphi_0^\ep||_{H^2(\ff{T})}\leq ||\varphi_0||_{H^2(\ff{T})}.$$
The above estimates provide that
$$||f^\ep_{x}||_{H^2(\ff{T})}\quad\text{is bounded}\quad \forall T<M(0)^{-1}\quad\text{uniformly in}\,\,\ep.$$
Using the Rellich's theorem we conclude the proof  taking the limit $\ep\rightarrow 0$.\quad\eop

\begin{remark}
We shall consider the equation
\begin{equation}
\pa_tf_{x}=-\pa_tg-ff_{xx}+(f_x)^2+gf_x - \nu\Lambda^{\alpha} f_x,\label{univis}
\end{equation}
which is (\ref{ecuacion}) with an extra dissipating term. For this system we have a local existence result  similar to theorem \ref{teorema1}.    Moreover, we can construct solutions that blow-up in finite time for $\alpha=1,\,2$ (see below) which show the existence of singularities for DPM with infinite energy.\quad\eop
\end{remark}

\smallskip

\subsection{Blow up.}

Next, we show that there exist a particular solution of
the equation (\ref{univis}), with $\nu\geq 0$ and $\alpha=1,\,2$, which blows up in finite time.

We consider  the following ansatz of (\ref{ecuacion}) and (\ref{univis})
\begin{equation}\label{fx}
f_{x}(x,t)=r(t)\cos(x),
\end{equation}
then $r$ satisfies
\begin{equation}\label{r}
\frac{dr(t)}{dt}=r(t)\int_0^tr^2(\tau)d\tau+\nu r.
\end{equation}
We define the
function
\begin{equation}\label{al}
\beta(t)=\int_0^tr^2(\tau)d\tau.
\end{equation}
Multiplying  the equation (\ref{r}) by $r(t)$
we have that $\beta$ satisfies
$$
\frac{d^2\beta(t)}{dt^2}=2\beta(t)\frac{d\beta(t)}{dt}+2\nu\frac{d\beta}{dt}.
$$
Integrating with respect to the variable $t$ yields
$$
\beta'(t)-\beta'(0)=\beta^2(t)-\beta^2(0)+2\nu(\beta-\beta(0)).
$$
Since  $\beta(0)=0$
and $\beta'(0)=r^2(0)$ we obtain
$$
 \frac{\beta'(t)}{r(0)^2+2\nu \beta + \beta(t)^2}=1.
$$
If we choose $r(0)^2>\nu^2$ it follows
\begin{equation*}
\beta(t)=\sqrt{r(0)^2-\nu^2}\tan\big(\sqrt{r(0)^2-\nu^2}\,t+\arctan(\frac{\nu}{\sqrt{r(0)^2-\nu^2}})\big)-\nu.
\end{equation*}
 Therefore, the function
\begin{equation}
f_{x}(x,t)=r(t) \cos(x),
\end{equation}
is a solution of equation (\ref{ecuacion}) and (\ref{univis}) which blows up at time
$$t=\frac{\big(\frac{\pi}{2}-\arctan(\frac{\nu}{\sqrt{r(0)^2-\nu^2}})\big)}{\sqrt{r(0)^2-\nu^2}}.$$

\smallskip\

\end{document}